\documentclass{amsart}
\usepackage{amsfonts}

\newtheorem{theorem}{Theorem}
\theoremstyle{plain}

\newtheorem{corollary}{Corollary}

\newtheorem{remark}{Remark}

\numberwithin{equation}{section}
\input{tcilatex}

\begin{document}
\title[Reverses of the Cauchy-Bunyakovsky-Schwarz Inequality]{Reverses of
the Cauchy-Bunyakovsky-Schwarz Inequality for $n-$tuples of Complex Numbers}
\author{S.S. Dragomir}
\address{School of Computer Science and Mathematics\\
Victoria University of Technology\\
PO Box 14428, MCMC 8001\\
Victoria, Australia.}
\email{sever.dragomir@vu.edu.au}
\urladdr{http://rgmia.vu.edu.au/SSDragomirWeb.html}
\date{November 07, 2003}

\begin{abstract}
Some new reverses of the Cauchy-Bunyakovsky-Schwarz inequality for $n-$%
tuples of real and complex numbers related to Cassels and Shisha-Mond
results are given.
\end{abstract}

\maketitle

\section{Introduction}

Let $\mathbf{\bar{a}}=\left( a_{1},\dots ,a_{n}\right) $ and $\mathbf{\bar{b}%
}=\left( b_{1},\dots ,b_{n}\right) $ be two positive $n$-tuples with the
property that there exists the positive numbers $m_{i},M_{i}$ $\left(
i=1,2\right) $ such that%
\begin{equation}
0<m_{1}\leq a_{i}\leq M_{1}<\infty \text{ \ and \ }0<m_{2}\leq b_{i}\leq
M_{2}<\infty ,  \label{1.1}
\end{equation}%
for each $i\in \left\{ 1,\dots ,n\right\} .$

The following reverses of the Cauchy-Bunyakovsky-Schwarz (CBS) inequality
are well known in the literature:

\begin{enumerate}
\item \textbf{P\'{o}lya-Szeg\"{o}'s inequality} \cite{8b}%
\begin{equation}
\frac{\sum_{k=1}^{n}a_{k}^{2}\sum_{k=1}^{n}b_{k}^{2}}{\left(
\sum_{k=1}^{n}a_{k}b_{k}\right) ^{2}}\leq \frac{1}{4}\left( \sqrt{\frac{%
M_{1}M_{2}}{m_{1}m_{2}}}+\sqrt{\frac{m_{1}m_{2}}{M_{1}M_{2}}}\right) ^{2};
\label{1.2}
\end{equation}

\item \textbf{Shisha-Mond's inequality} \cite{9b}%
\begin{equation}
\frac{\sum_{k=1}^{n}a_{k}^{2}}{\sum_{k=1}^{n}a_{k}b_{k}}-\frac{%
\sum_{k=1}^{n}a_{k}b_{k}}{\sum_{k=1}^{n}b_{k}^{2}}\leq \left( \sqrt{\frac{%
M_{1}}{m_{2}}}-\sqrt{\frac{m_{1}}{M_{2}}}\right) ^{2};  \label{1.3}
\end{equation}

\item \textbf{Ozeki's inequality} \cite{7b}%
\begin{equation}
\sum_{k=1}^{n}a_{k}^{2}\sum_{k=1}^{n}b_{k}^{2}-\left(
\sum_{k=1}^{n}a_{k}b_{k}\right) ^{2}\leq \frac{1}{4}n^{2}\left(
M_{1}M_{2}-m_{1}m_{2}\right) ^{2};  \label{1.3a}
\end{equation}

\item \textbf{Diaz-Metcalf's inequality} \cite{1b}%
\begin{equation}
\sum_{k=1}^{n}b_{k}^{2}+\frac{m_{2}M_{2}}{m_{1}M_{1}}\sum_{k=1}^{n}a_{k}^{2}%
\leq \left( \frac{M_{2}}{m_{1}}+\frac{m_{2}}{M_{1}}\right)
\sum_{k=1}^{n}a_{k}b_{k}.  \label{1.4}
\end{equation}
\end{enumerate}

If the weight $\mathbf{\bar{w}}=\left( w_{1},\dots ,w_{n}\right) $ is a
positive $n-$tuple, then we have the following inequalities, which are also
well known.

\begin{enumerate}
\item[5.] \textbf{Cassel's inequality} \cite{10b}\newline
If the positive $n-$tuples $\mathbf{\bar{a}}=\left( a_{1},\dots
,a_{n}\right) $ and $\mathbf{\bar{b}}=\left( b_{1},\dots ,b_{n}\right) $
satisfy the condition%
\begin{equation}
0<m\leq \frac{a_{k}}{b_{k}}\leq M<\infty \text{ \ for each \ }k\in \left\{
1,\dots ,n\right\} ,  \label{1.5}
\end{equation}%
where $m,M$ are given, then%
\begin{equation}
\frac{\sum_{k=1}^{n}w_{k}a_{k}^{2}\sum_{k=1}^{n}w_{k}b_{k}^{2}}{\left(
\sum_{k=1}^{n}w_{k}a_{k}b_{k}\right) ^{2}}\leq \frac{\left( M+m\right) ^{2}}{%
4mM}.  \label{1.6}
\end{equation}

\item[6.] \textbf{Grueb-Reinboldt's inequality} \cite{4b}\newline
If $\mathbf{\bar{a}}$ and $\mathbf{\bar{b}}$ satisfy the condition (\ref{1.1}%
), then%
\begin{equation}
\frac{\sum_{k=1}^{n}w_{k}a_{k}^{2}\sum_{k=1}^{n}w_{k}b_{k}^{2}}{\left(
\sum_{k=1}^{n}w_{k}a_{k}b_{k}\right) ^{2}}\leq \frac{\left(
M_{1}M_{2}+m_{1}m_{2}\right) ^{2}}{4m_{1}m_{2}M_{1}M_{2}}.  \label{1.7}
\end{equation}

\item[7.] \textbf{Generalised Diaz-Metcalf inequality} \cite{1b} (see also 
\cite[p. 123]{6b})\newline
If $u,v\in \left[ 0,1\right] $ and $v\leq u,$ $u+v=1$ and (\ref{1.5}) holds,
then one has the inequality%
\begin{equation}
u\sum_{k=1}^{n}w_{k}b_{k}^{2}+vmM\sum_{k=1}^{n}w_{k}a_{k}^{2}\leq \left(
vm+uM\right) \sum_{k=1}^{n}w_{k}a_{k}b_{k}.  \label{1.8}
\end{equation}

\item[8.] \textbf{Klamkin-McLenaghan's inequality} \cite{5b}

If $\mathbf{\bar{a}}$ and $\mathbf{\bar{b}}$ satisfy (\ref{1.5}), then we
have the inequality%
\begin{equation}
\sum_{k=1}^{n}w_{k}a_{k}^{2}\sum_{k=1}^{n}w_{k}b_{k}^{2}-\left(
\sum_{k=1}^{n}w_{k}a_{k}b_{k}\right) ^{2}\leq \left( \sqrt{M}-\sqrt{m}%
\right) ^{2}\sum_{k=1}^{n}w_{k}a_{k}b_{k}\sum_{k=1}^{n}w_{k}a_{k}^{2}.
\label{1.9}
\end{equation}
\end{enumerate}

For other reverse results of the (CBS)-inequality, see the recent survey
online \cite{3b}.

The main aim of this paper is to point out some new reverse inequalities of
the classical Cauchy-Bunyakovsky-Schwarz result for both real and complex $%
n- $tuples.

\section{Some Reverses of the Cauchy-Bunyakovsky-Schwarz Inequality}

The following result holds.

\begin{theorem}
\label{t2.1}Let $\mathbf{\bar{a}}=\left( a_{1},\dots ,a_{n}\right) $, $%
\mathbf{\bar{b}}=\left( b_{1},\dots ,b_{n}\right) \in \mathbb{K}^{n},$ where 
$\mathbb{K}=\mathbb{R},\mathbb{C}$ and $\mathbf{\bar{p}}=\left( p_{1},\dots
,p_{n}\right) \in \mathbb{R}_{+}^{n}$ with $\sum_{i=1}^{n}p_{i}=1.$ If $%
b_{i}\neq 0,$ $i\in \left\{ 1,\dots ,n\right\} $ and there exists the
constant $\alpha \in \mathbb{K}$ and $r>0$ such that for any $k\in \left\{
1,\dots ,n\right\} $%
\begin{equation}
\frac{a_{k}}{\overline{b_{k}}}\in \bar{D}\left( \alpha ,r\right) :=\left\{
z\in \mathbb{K}|\left\vert z-\alpha \right\vert \leq r\right\} ,  \label{2.1}
\end{equation}%
then we have the inequality%
\begin{equation}
\sum_{k=1}^{n}p_{k}\left\vert a_{k}\right\vert ^{2}+\left( \left\vert \alpha
\right\vert ^{2}-r^{2}\right) \sum_{k=1}^{n}p_{k}\left\vert b_{k}\right\vert
^{2}\leq 2\func{Re}\left[ \bar{\alpha}\left(
\sum_{k=1}^{n}p_{k}a_{k}b_{k}\right) \right] .  \label{2.2}
\end{equation}%
The constant $c=2$ is best possible in the sense that it cannot be replaced
by a smaller constant.
\end{theorem}

\begin{proof}
From (\ref{2.1}) we have $\left\vert a_{k}-\alpha \bar{b}_{k}\right\vert
^{2}\leq r\left\vert b_{k}\right\vert ^{2}$ for each $k\in \left\{ 1,\dots
,n\right\} ,$ which is clearly equivalent to%
\begin{equation}
\left\vert a_{k}\right\vert ^{2}+\left( \left\vert \alpha \right\vert
^{2}-r^{2}\right) \left\vert b_{k}\right\vert ^{2}\leq 2\func{Re}\left[ \bar{%
\alpha}\left( a_{k}b_{k}\right) \right]  \label{2.3}
\end{equation}%
for each $k\in \left\{ 1,\dots ,n\right\} .$

Multiplying (\ref{2.3}) with $p_{k}\geq 0$ \ and summing over $k$ from $1$
to $n,$ we deduce the first inequality in (\ref{1.2}). The second inequality
is obvious.

To prove the sharpness of the constant $2,$ assume that under the hypothesis
of the theorem there exists a constant $c>0$ such that%
\begin{equation}
\sum_{k=1}^{n}p_{k}\left\vert a_{k}\right\vert ^{2}+\left( \left\vert \alpha
\right\vert ^{2}-r^{2}\right) \sum_{k=1}^{n}p_{k}\left\vert b_{k}\right\vert
^{2}\leq c\func{Re}\left[ \bar{\alpha}\left(
\sum_{k=1}^{n}p_{k}a_{k}b_{k}\right) \right] ,  \label{2.4}
\end{equation}%
provided $\frac{a_{k}}{\overline{b_{k}}}\in \bar{D}\left( \alpha ,r\right) ,$
$k\in \left\{ 1,\dots ,n\right\} .$

Assume that $n=2,$ $p_{1}=p_{2}=\frac{1}{2},$ $b_{1}=b_{2}=1,$ $\alpha =r>0$
and $a_{2}=2r,$ $a_{1}=0.$ Then $\left\vert \frac{a_{2}}{b_{2}}-\alpha
\right\vert =r,$ $\left\vert \frac{a_{1}}{b_{1}}-\alpha \right\vert =r$
showing that the condition (\ref{2.1}) holds. For these choices, the
inequality (\ref{2.4}) becomes $2r^{2}\leq cr^{2},$ giving $c\geq 2.$
\end{proof}

The case where the disk $\bar{D}\left( \alpha ,r\right) $ does not contain
the origin, i.e., $\left\vert \alpha \right\vert >r,$ provides the following
interesting reverse of the Cauchy-Bunyakovsky-Schwarz inequality.

\begin{theorem}
\label{t2.2}Let $\mathbf{\bar{a}}$, $\mathbf{\bar{b}}$, $\mathbf{\bar{p}}$
as in Theorem \ref{t2.1} and assume that $\left\vert \alpha \right\vert
>r>0. $ Then we have the inequality%
\begin{align}
\sum_{k=1}^{n}p_{k}\left\vert a_{k}\right\vert
^{2}\sum_{k=1}^{n}p_{k}\left\vert b_{k}\right\vert ^{2}& \leq \frac{1}{%
\left\vert \alpha \right\vert ^{2}-r^{2}}\left\{ \func{Re}\left[ \bar{\alpha}%
\left( \sum_{k=1}^{n}p_{k}a_{k}b_{k}\right) \right] \right\} ^{2}
\label{2.5} \\
& \leq \frac{\left\vert \alpha \right\vert ^{2}}{\left\vert \alpha
\right\vert ^{2}-r^{2}}\left\vert \sum_{k=1}^{n}p_{k}a_{k}b_{k}\right\vert
^{2}.  \notag
\end{align}%
The constant $c=1$ in the first and second inequality is best possible in
the sense that it cannot be replaced by a smaller constant.
\end{theorem}

\begin{proof}
Since $\left\vert \alpha \right\vert >r,$ we may divide (\ref{2.2}) by $%
\sqrt{\left\vert \alpha \right\vert ^{2}-r^{2}}>0$ to obtain%
\begin{multline}
\frac{1}{\sqrt{\left\vert \alpha \right\vert ^{2}-r^{2}}}\sum_{k=1}^{n}p_{k}%
\left\vert a_{k}\right\vert ^{2}+\sqrt{\left\vert \alpha \right\vert
^{2}-r^{2}}\sum_{k=1}^{n}p_{k}\left\vert b_{k}\right\vert ^{2}  \label{2.6}
\\
\leq \frac{2}{\sqrt{\left\vert \alpha \right\vert ^{2}-r^{2}}}\func{Re}\left[
\bar{\alpha}\left( \sum_{k=1}^{n}p_{k}a_{k}b_{k}\right) \right] .
\end{multline}%
On the other hand, by the use of the following elementary inequality%
\begin{equation}
\frac{1}{\beta }p+\beta q\geq 2\sqrt{pq}\text{ \ for \ }\beta >0\text{ \ and
\ }p,q\geq 0,  \label{2.6.a}
\end{equation}%
we may state that%
\begin{multline}
2\left( \sum_{k=1}^{n}p_{k}\left\vert a_{k}\right\vert ^{2}\right) ^{\frac{1%
}{2}}\cdot \left( \sum_{k=1}^{n}p_{k}\left\vert b_{k}\right\vert ^{2}\right)
^{\frac{1}{2}}  \label{2.7} \\
\leq \frac{1}{\sqrt{\left\vert \alpha \right\vert ^{2}-r^{2}}}%
\sum_{k=1}^{n}p_{k}\left\vert a_{k}\right\vert ^{2}+\sqrt{\left\vert \alpha
\right\vert ^{2}-r^{2}}\sum_{k=1}^{n}p_{k}\left\vert b_{k}\right\vert ^{2}.
\end{multline}%
Utilising (\ref{2.6}) and (\ref{2.7}), we deduce%
\begin{equation*}
\left( \sum_{k=1}^{n}p_{k}\left\vert a_{k}\right\vert ^{2}\right) ^{\frac{1}{%
2}}\cdot \left( \sum_{k=1}^{n}p_{k}\left\vert b_{k}\right\vert ^{2}\right) ^{%
\frac{1}{2}}\leq \frac{1}{\sqrt{\left\vert \alpha \right\vert ^{2}-r^{2}}}%
\func{Re}\left[ \bar{\alpha}\left( \sum_{k=1}^{n}p_{k}a_{k}b_{k}\right) %
\right] ,
\end{equation*}%
which is clearly equivalent to the first inequality in (\ref{2.6}).

The second inequality is obvious.

To prove the sharpness of the constant, assume that (\ref{2.5}) holds with a
constant $c>0,$ i.e.,%
\begin{equation}
\sum_{k=1}^{n}p_{k}\left\vert a_{k}\right\vert
^{2}\sum_{k=1}^{n}p_{k}\left\vert b_{k}\right\vert ^{2}\leq \frac{c}{%
\left\vert \alpha \right\vert ^{2}-r^{2}}\left\{ \func{Re}\left[ \bar{\alpha}%
\left( \sum_{k=1}^{n}p_{k}a_{k}b_{k}\right) \right] \right\} ^{2}
\label{2.8}
\end{equation}%
provided $\frac{a_{k}}{\overline{b_{k}}}\in \bar{D}\left( \alpha ,r\right) $
and $\left\vert \alpha \right\vert >1.$

For $n=2,$ $b_{2}=b_{1}=1,$ $p_{1}=p_{2}=\frac{1}{2},$ $a_{2},a_{1}\in 
\mathbb{R}$, $\alpha ,r>0$ and $\alpha >r,$ we get from (\ref{2.8}) that%
\begin{equation}
\frac{a_{1}^{2}+a_{2}^{2}}{2}\leq \frac{c\alpha ^{2}}{\alpha ^{2}-r^{2}}%
\left( \frac{a_{1}+a_{2}}{2}\right) ^{2}.  \label{2.9}
\end{equation}%
If we choose $a_{2}=\alpha +r,$ $a_{1}=\alpha -r,$ then $\left\vert
a_{i}-\alpha \right\vert \leq r,$ $i=1,2$ and by (\ref{2.9}) we deduce%
\begin{equation*}
\alpha ^{2}+r^{2}\leq \frac{c\alpha ^{4}}{\alpha ^{2}-r^{2}},
\end{equation*}%
which is clearly equivalent to%
\begin{equation*}
\left( c-1\right) \alpha ^{4}+r^{4}\geq 0\text{ \ for \ }\alpha >r>0.
\end{equation*}%
If in this inequality we choose $\alpha =1,$ $r=\varepsilon \in \left(
0,1\right) $ and let $\varepsilon \rightarrow 0+,$ then we deduce $c\geq 1.$
\end{proof}

The following corollary is a natural consequence of the above theorem.

\begin{corollary}
\label{c2.3}Under the assumptions of Theorem \ref{t2.2}, we have the
following additive reverse of the Cauchy-Bunyakovsky-Schwarz inequality%
\begin{align}
0& \leq \sum_{k=1}^{n}p_{k}\left\vert a_{k}\right\vert
^{2}\sum_{k=1}^{n}p_{k}\left\vert b_{k}\right\vert ^{2}-\left\vert
\sum_{k=1}^{n}p_{k}a_{k}b_{k}\right\vert ^{2}  \label{2.10} \\
& \leq \frac{r^{2}}{\left\vert \alpha \right\vert ^{2}-r^{2}}\left\vert
\sum_{k=1}^{n}p_{k}a_{k}b_{k}\right\vert ^{2}.  \notag
\end{align}%
The constant $c=1$ is best possible in the sense mentioned above.
\end{corollary}

\begin{remark}
\label{r2.4}If in Theorem \ref{t2.1}, we assume that $\left\vert \alpha
\right\vert =r,$ then we obtain the inequality:%
\begin{align}
\sum_{k=1}^{n}p_{k}\left\vert a_{k}\right\vert ^{2}& \leq 2\func{Re}\left[ 
\bar{\alpha}\left( \sum_{k=1}^{n}p_{k}a_{k}b_{k}\right) \right]  \label{2.11}
\\
& \leq 2\left\vert \alpha \right\vert \left\vert
\sum_{k=1}^{n}p_{k}a_{k}b_{k}\right\vert .  \notag
\end{align}%
The constant $2$ is sharp in both inequalities.

We also remark that, if $r>\left\vert \alpha \right\vert $, then (\ref{2.2})
may be written as%
\begin{align}
\sum_{k=1}^{n}p_{k}\left\vert a_{k}\right\vert ^{2}& \leq \left(
r^{2}-\left\vert \alpha \right\vert ^{2}\right)
\sum_{k=1}^{n}p_{k}\left\vert b_{k}\right\vert ^{2}+2\func{Re}\left[ \bar{%
\alpha}\left( \sum_{k=1}^{n}p_{k}a_{k}b_{k}\right) \right]  \label{2.12} \\
& \leq \left( r^{2}-\left\vert \alpha \right\vert ^{2}\right)
\sum_{k=1}^{n}p_{k}\left\vert b_{k}\right\vert ^{2}+2\left\vert \alpha
\right\vert \left\vert \sum_{k=1}^{n}p_{k}a_{k}b_{k}\right\vert .  \notag
\end{align}
\end{remark}

The following reverse of the Cauchy-Bunyakovsky-Schwarz inequality also
holds.

\begin{theorem}
\label{t2.4}Let $\mathbf{\bar{a}}$, $\mathbf{\bar{b}}$, $\mathbf{\bar{p}}$
be as in Theorem \ref{t2.1} and assume that $\alpha \in \mathbb{K}$, $\alpha
\neq 0$ and $r>0.$ Then we have the inequalities%
\begin{align}
0& \leq \left( \sum_{k=1}^{n}p_{k}\left\vert a_{k}\right\vert ^{2}\right) ^{%
\frac{1}{2}}\cdot \left( \sum_{k=1}^{n}p_{k}\left\vert b_{k}\right\vert
^{2}\right) ^{\frac{1}{2}}-\left\vert
\sum_{k=1}^{n}p_{k}a_{k}b_{k}\right\vert   \label{2.11a} \\
& \leq \left( \sum_{k=1}^{n}p_{k}\left\vert a_{k}\right\vert ^{2}\right) ^{%
\frac{1}{2}}\cdot \left( \sum_{k=1}^{n}p_{k}\left\vert b_{k}\right\vert
^{2}\right) ^{\frac{1}{2}}-\func{Re}\left[ \frac{\bar{\alpha}}{\left\vert
\alpha \right\vert }\left( \sum_{k=1}^{n}p_{k}a_{k}b_{k}\right) \right]  
\notag \\
& \leq \frac{1}{2}\cdot \frac{r^{2}}{\left\vert \alpha \right\vert }%
\sum_{k=1}^{n}p_{k}\left\vert b_{k}\right\vert ^{2}.  \notag
\end{align}%
The constant $\frac{1}{2}$ is best possible in the sense mentioned above.
\end{theorem}

\begin{proof}
From Theorem \ref{t2.1}, we have%
\begin{equation}
\sum_{k=1}^{n}p_{k}\left\vert a_{k}\right\vert ^{2}+\left\vert \alpha
\right\vert ^{2}\sum_{k=1}^{n}p_{k}\left\vert b_{k}\right\vert ^{2}\leq 2%
\func{Re}\left[ \bar{\alpha}\left( \sum_{k=1}^{n}p_{k}a_{k}b_{k}\right) %
\right] +r^{2}\sum_{k=1}^{n}p_{k}\left\vert b_{k}\right\vert ^{2}.
\label{2.12a}
\end{equation}%
Since $\alpha \neq 0,$ we can divide (\ref{2.12a}) by $\left\vert \alpha
\right\vert ,$ getting%
\begin{multline}
\frac{1}{\left\vert \alpha \right\vert }\sum_{k=1}^{n}p_{k}\left\vert
a_{k}\right\vert ^{2}+\left\vert \alpha \right\vert
\sum_{k=1}^{n}p_{k}\left\vert b_{k}\right\vert ^{2}  \label{2.13} \\
\leq 2\func{Re}\left[ \frac{\bar{\alpha}}{\left\vert \alpha \right\vert }%
\left( \sum_{k=1}^{n}p_{k}a_{k}b_{k}\right) \right] +\frac{r^{2}}{\left\vert
\alpha \right\vert }\sum_{k=1}^{n}p_{k}\left\vert b_{k}\right\vert ^{2}.
\end{multline}%
Utilising the inequality (\ref{2.6.a}), we may state that%
\begin{equation}
2\left( \sum_{k=1}^{n}p_{k}\left\vert a_{k}\right\vert ^{2}\right) ^{\frac{1%
}{2}}\cdot \left( \sum_{k=1}^{n}p_{k}\left\vert b_{k}\right\vert ^{2}\right)
^{\frac{1}{2}}\leq \frac{1}{\left\vert \alpha \right\vert }%
\sum_{k=1}^{n}p_{k}\left\vert a_{k}\right\vert ^{2}+\left\vert \alpha
\right\vert \sum_{k=1}^{n}p_{k}\left\vert b_{k}\right\vert ^{2}.
\label{2.14}
\end{equation}%
Making use of (\ref{2.13}) and (\ref{2.14}), we deduce the second inequality
in (\ref{2.11a}).

The first inequality in (\ref{2.11a}) is obvious.

To prove the sharpness of the constant $\frac{1}{2},$ assume that there
exists a $c>0$ such that%
\begin{multline}
\left( \sum_{k=1}^{n}p_{k}\left\vert a_{k}\right\vert ^{2}\right) ^{\frac{1}{%
2}}\cdot \left( \sum_{k=1}^{n}p_{k}\left\vert b_{k}\right\vert ^{2}\right) ^{%
\frac{1}{2}}-\func{Re}\left[ \frac{\bar{\alpha}}{\left\vert \alpha
\right\vert }\left( \sum_{k=1}^{n}p_{k}a_{k}b_{k}\right) \right]
\label{2.15} \\
\leq c\cdot \frac{r^{2}}{\left\vert \alpha \right\vert }\sum_{k=1}^{n}p_{k}%
\left\vert b_{k}\right\vert ^{2},
\end{multline}%
provided $\left\vert \frac{a_{k}}{\overline{b_{k}}}-\alpha \right\vert \leq
r,$ $\alpha \neq 0,$ $r>0.$

If we choose $n=2,$ $\alpha >0,$ $b_{1}=b_{2}=1,$ $a_{1}=\alpha +r,$ $%
a_{2}=\alpha -r,$ then from (\ref{2.15}) we deduce%
\begin{equation}
\sqrt{r^{2}+\alpha ^{2}}-\alpha \leq c\frac{r^{2}}{\alpha }.  \label{2.16}
\end{equation}%
If we multiply (\ref{2.16}) with $\sqrt{r^{2}+\alpha ^{2}}+\alpha >0$ and
then divide it by $r>0,$ we deduce%
\begin{equation}
1\leq \frac{\sqrt{r^{2}+\alpha ^{2}}+\alpha }{\alpha }\cdot c  \label{2.17}
\end{equation}%
for any $r>0,$ $\alpha >0.$

If in (\ref{2.17}) we let $r\rightarrow 0+,$ then we get $c\geq \frac{1}{2},$
and the sharpness of the constant is proved.
\end{proof}

\section{A Cassels Type Inequality for Complex Numbers}

The following result holds.

\begin{theorem}
\label{t3.1}Let $\mathbf{\bar{a}}=\left( a_{1},\dots ,a_{n}\right) $, $%
\mathbf{\bar{b}}=\left( b_{1},\dots ,b_{n}\right) \in \mathbb{K}^{n},$ where 
$\mathbb{K}=\mathbb{R},\mathbb{C}$ and $\mathbf{\bar{p}}=\left( p_{1},\dots
,p_{n}\right) \in \mathbb{R}_{+}^{n}$ with $\sum_{i=1}^{n}p_{i}=1.$ If $%
b_{i}\neq 0,$ $i\in \left\{ 1,\dots ,n\right\} $ and there exist the
constants $\gamma ,\Gamma \in \mathbb{K}$ with $\func{Re}\left( \Gamma \bar{%
\gamma}\right) >0$ and $\Gamma \neq \gamma ,$ so that either%
\begin{equation}
\left\vert \frac{a_{k}}{\overline{b_{k}}}-\frac{\gamma +\Gamma }{2}%
\right\vert \leq \frac{1}{2}\left\vert \Gamma -\gamma \right\vert \text{ \
for each \ }k\in \left\{ 1,\dots ,n\right\} ,  \label{3.1}
\end{equation}%
or, equivalently,%
\begin{equation}
\func{Re}\left[ \left( \Gamma -\frac{a_{k}}{b_{k}}\right) \left( \frac{%
\overline{a_{k}}}{b_{k}}-\bar{\gamma}\right) \right] \geq 0\text{ \ for each
\ }k\in \left\{ 1,\dots ,n\right\}   \label{3.2}
\end{equation}%
holds, then we have the inequalities%
\begin{align}
\sum_{k=1}^{n}p_{k}\left\vert a_{k}\right\vert
^{2}\sum_{k=1}^{n}p_{k}\left\vert b_{k}\right\vert ^{2}& \leq \frac{1}{2%
\func{Re}\left( \Gamma \bar{\gamma}\right) }\left\{ \func{Re}\left[ \left( 
\bar{\gamma}+\bar{\Gamma}\right) \sum_{k=1}^{n}p_{k}a_{k}b_{k}\right]
\right\} ^{2}  \label{3.3} \\
& \leq \frac{\left\vert \Gamma +\gamma \right\vert ^{2}}{4\func{Re}\left(
\Gamma \bar{\gamma}\right) }\left\vert
\sum_{k=1}^{n}p_{k}a_{k}b_{k}\right\vert ^{2}.  \notag
\end{align}%
The constants $\frac{1}{2}$ and $\frac{1}{4}$ are best possible in (\ref{3.3}%
).
\end{theorem}

\begin{proof}
The fact that the relations (\ref{3.1}) and (\ref{3.2}) are equivalent
follows by the simple fact that for $z,u,U\in \mathbb{C}$, the following
inequalities are equivalent%
\begin{equation*}
\left\vert z-\frac{u+U}{2}\right\vert \leq \frac{1}{2}\left\vert
U-u\right\vert
\end{equation*}%
and%
\begin{equation*}
\func{Re}\left[ \left( u-z\right) \left( \bar{z}-\bar{u}\right) \right] \geq
0.
\end{equation*}%
Define $\alpha =\frac{\gamma +\Gamma }{2}$ and $r=\frac{1}{2}\left\vert
\Gamma -\gamma \right\vert .$ Then%
\begin{equation*}
\left\vert \alpha \right\vert ^{2}-r^{2}=\frac{\left\vert \Gamma +\gamma
\right\vert ^{2}}{4}-\frac{\left\vert \Gamma -\gamma \right\vert ^{2}}{4}=%
\func{Re}\left( \Gamma \bar{\gamma}\right) >0.
\end{equation*}%
Consequently, we may apply Theorem \ref{t2.2}, and the inequalities (\ref%
{3.3}) are proved.

The sharpness of the constants may be proven in a similar way to that in the
proof of Theorem \ref{t2.2}, and we omit the details.
\end{proof}

The following additive version also holds.

\begin{corollary}
\label{c3.2}With the assumptions in Theorem \ref{t3.1}, we have%
\begin{equation}
\sum_{k=1}^{n}p_{k}\left\vert a_{k}\right\vert
^{2}\sum_{k=1}^{n}p_{k}\left\vert b_{k}\right\vert ^{2}-\left\vert
\sum_{k=1}^{n}p_{k}a_{k}b_{k}\right\vert ^{2}\leq \frac{\left\vert \Gamma
-\gamma \right\vert ^{2}}{4\func{Re}\left( \Gamma \bar{\gamma}\right) }%
\left\vert \sum_{k=1}^{n}p_{k}a_{k}b_{k}\right\vert ^{2}.  \label{3.4}
\end{equation}%
The constant $\frac{1}{4}$ is also best possible.
\end{corollary}

\begin{remark}
\label{r3.3}With the above assumptions and if $\func{Re}\left( \Gamma \bar{%
\gamma}\right) =0,$ then by the use of Remark \ref{r2.4}, we may deduce the
inequality%
\begin{equation}
\sum_{k=1}^{n}p_{k}\left\vert a_{k}\right\vert ^{2}\leq \func{Re}\left[
\left( \bar{\gamma}+\bar{\Gamma}\right) \sum_{k=1}^{n}p_{k}a_{k}b_{k}\right]
\leq \left\vert \Gamma +\gamma \right\vert \left\vert
\sum_{k=1}^{n}p_{k}a_{k}b_{k}\right\vert .  \label{3.5}
\end{equation}%
If $\func{Re}\left( \Gamma \bar{\gamma}\right) <0,$ then, by Remark \ref%
{r2.4}, we also have%
\begin{align}
\sum_{k=1}^{n}p_{k}\left\vert a_{k}\right\vert ^{2}& \leq -\func{Re}\left(
\Gamma \bar{\gamma}\right) \sum_{k=1}^{n}p_{k}\left\vert b_{k}\right\vert
^{2}+\func{Re}\left[ \left( \bar{\Gamma}+\bar{\gamma}\right)
\sum_{k=1}^{n}p_{k}a_{k}b_{k}\right]  \label{3.6} \\
& \leq -\func{Re}\left( \Gamma \bar{\gamma}\right)
\sum_{k=1}^{n}p_{k}\left\vert b_{k}\right\vert ^{2}+\left\vert \Gamma
+\gamma \right\vert \left\vert \sum_{k=1}^{n}p_{k}a_{k}b_{k}\right\vert . 
\notag
\end{align}
\end{remark}

\begin{remark}
\label{r3.4}If $a_{k},b_{k}>0$ and there exist the constants $m,M>0$ $\left(
M>m\right) $ with%
\begin{equation}
m\leq \frac{a_{k}}{b_{k}}\leq M\text{ \ for each \ }k\in \left\{ 1,\dots
,n\right\} ,  \label{3.7}
\end{equation}%
then, obviously (\ref{3.1}) holds with $\gamma =m,$ $\Gamma =M,$ also $%
\Gamma \bar{\gamma}=Mm>0$ and by (\ref{3.3}) we deduce%
\begin{equation}
\sum_{k=1}^{n}p_{k}a_{k}^{2}\sum_{k=1}^{n}p_{k}b_{k}^{2}\leq \frac{\left(
M+m\right) ^{2}}{4mM}\left( \sum_{k=1}^{n}p_{k}a_{k}b_{k}\right) ^{2},
\label{3.8}
\end{equation}%
that is, Cassels' inequality.
\end{remark}

\section{A Shisha-Mond Type Inequality for Complex Numbers}

The following result holds.

\begin{theorem}
\label{t4.1}Let $\mathbf{\bar{a}}=\left( a_{1},\dots ,a_{n}\right) $, $%
\mathbf{\bar{b}}=\left( b_{1},\dots ,b_{n}\right) \in \mathbb{K}^{n},$ where 
$\mathbb{K}=\mathbb{R},\mathbb{C}$ and $\mathbf{\bar{p}}=\left( p_{1},\dots
,p_{n}\right) \in \mathbb{R}_{+}^{n}$ with $\sum_{i=1}^{n}p_{i}=1.$ If $%
b_{i}\neq 0,$ $i\in \left\{ 1,\dots ,n\right\} $ and there exist the
constants $\gamma ,\Gamma \in \mathbb{K}$ such that $\Gamma \neq \gamma
,-\gamma $ and either%
\begin{equation}
\left\vert \frac{a_{k}}{\overline{b_{k}}}-\frac{\gamma +\Gamma }{2}%
\right\vert \leq \frac{1}{2}\left\vert \Gamma -\gamma \right\vert \text{ \
for each \ }k\in \left\{ 1,\dots ,n\right\} ,  \label{4.1}
\end{equation}%
or, equivalently,%
\begin{equation}
\func{Re}\left[ \left( \Gamma -\frac{a_{k}}{b_{k}}\right) \left( \frac{%
\overline{a_{k}}}{b_{k}}-\bar{\gamma}\right) \right] \geq 0\text{ \ for each
\ }k\in \left\{ 1,\dots ,n\right\} ,  \label{4.2}
\end{equation}%
holds, then we have the inequalities%
\begin{align}
0& \leq \left( \sum_{k=1}^{n}p_{k}\left\vert a_{k}\right\vert ^{2}\right) ^{%
\frac{1}{2}}\cdot \left( \sum_{k=1}^{n}p_{k}\left\vert b_{k}\right\vert
^{2}\right) ^{\frac{1}{2}}-\left\vert
\sum_{k=1}^{n}p_{k}a_{k}b_{k}\right\vert  \label{4.3} \\
& \leq \left( \sum_{k=1}^{n}p_{k}\left\vert a_{k}\right\vert ^{2}\right) ^{%
\frac{1}{2}}\cdot \left( \sum_{k=1}^{n}p_{k}\left\vert b_{k}\right\vert
^{2}\right) ^{\frac{1}{2}}-\func{Re}\left[ \frac{\bar{\Gamma}+\bar{\gamma}}{%
\left\vert \Gamma +\gamma \right\vert }\sum_{k=1}^{n}p_{k}a_{k}b_{k}\right] 
\notag \\
& \leq \frac{1}{4}\cdot \frac{\left\vert \Gamma -\gamma \right\vert ^{2}}{%
\left\vert \Gamma +\gamma \right\vert }\sum_{k=1}^{n}p_{k}\left\vert
b_{k}\right\vert ^{2}.  \notag
\end{align}%
The constant $\frac{1}{4}$ is best possible in the sense that it cannot be
replaced by a smaller constant.
\end{theorem}

\begin{proof}
Follows by Theorem \ref{t2.4} on choosing $\alpha =\frac{\gamma +\Gamma }{2}%
\neq 0$ and $r=\frac{1}{2}\left\vert \Gamma -\gamma \right\vert >0.$

The proof for the best constant follows in a similar way to that in the
proof of Theorem \ref{t2.4} and we omit the details.
\end{proof}

\begin{remark}
\label{r4.2}If $a_{k},b_{k}>0$ and there exists the constants $m,M>0$ $%
\left( M>m\right) $ with%
\begin{equation}
m\leq \frac{a_{k}}{b_{k}}\leq M\text{ \ for each \ }k\in \left\{ 1,\dots
,n\right\} ,  \label{4.4}
\end{equation}%
then we have the inequality%
\begin{align}
0& \leq \left( \sum_{k=1}^{n}p_{k}a_{k}^{2}\right) ^{\frac{1}{2}}\cdot
\left( \sum_{k=1}^{n}p_{k}b_{k}^{2}\right) ^{\frac{1}{2}}-%
\sum_{k=1}^{n}p_{k}a_{k}b_{k}  \label{4.5} \\
& \leq \frac{1}{4}\cdot \frac{\left( M-m\right) ^{2}}{\left( M+m\right) }%
\sum_{k=1}^{n}p_{k}b_{k}^{2}.  \notag
\end{align}%
The constant $\frac{1}{4}$ is best possible. For $p_{k}=\frac{1}{n},$ $k\in
\left\{ 1,\dots ,n\right\} ,$ we recapture the result from \cite[Theorem 5.21%
]{2b} that has been obtained from a reverse inequality due to Shisha and
Mond \cite{8b}.
\end{remark}

\section{Further Reverses of the (CBS)-Inequality}

The following result holds.

\begin{theorem}
\label{t5.1}Let $\mathbf{\bar{a}}=\left( a_{1},\dots ,a_{n}\right) $, $%
\mathbf{\bar{b}}=\left( b_{1},\dots ,b_{n}\right) \in \mathbb{K}^{n}$ and $%
r>0$ such that for $p_{i}\geq 0$ with $\sum_{i=1}^{n}p_{i}=1$%
\begin{equation}
\sum_{i=1}^{n}p_{i}\left\vert b_{i}-\overline{a_{i}}\right\vert ^{2}\leq
r^{2}<\sum_{i=1}^{n}p_{i}\left\vert a_{i}\right\vert ^{2}.  \label{5.1}
\end{equation}%
Then we have the inequality%
\begin{align}
0& \leq \sum_{i=1}^{n}p_{i}\left\vert a_{i}\right\vert
^{2}\sum_{i=1}^{n}p_{i}\left\vert b_{i}\right\vert ^{2}-\left\vert
\sum_{i=1}^{n}p_{i}a_{i}b_{i}\right\vert ^{2}  \label{5.2} \\
& \leq \sum_{i=1}^{n}p_{i}\left\vert a_{i}\right\vert
^{2}\sum_{i=1}^{n}p_{i}\left\vert b_{i}\right\vert ^{2}-\left[ \func{Re}%
\left( \sum_{i=1}^{n}p_{i}a_{i}b_{i}\right) \right] ^{2}  \notag \\
& \leq r^{2}\sum_{i=1}^{n}p_{i}\left\vert b_{i}\right\vert ^{2}.  \notag
\end{align}%
The constant $c=1$ in front of $r^{2}$ is best possible in the sense that it
cannot be replaced by a smaller constant.
\end{theorem}

\begin{proof}
From the first condition in (\ref{5.1}), we have%
\begin{equation*}
\sum_{i=1}^{n}p_{i}\left[ \left\vert b_{i}\right\vert ^{2}-2\func{Re}\left(
b_{i}a_{i}\right) +\left\vert a_{i}\right\vert ^{2}\right] \leq r^{2},
\end{equation*}%
giving%
\begin{equation}
\sum_{i=1}^{n}p_{i}\left\vert b_{i}\right\vert
^{2}+\sum_{i=1}^{n}p_{i}\left\vert a_{i}\right\vert ^{2}-r^{2}\leq 2\func{Re}%
\left( \sum_{i=1}^{n}p_{i}a_{i}b_{i}\right) .  \label{5.3}
\end{equation}%
Since, by the second condition in (\ref{5.1}) we have%
\begin{equation*}
\sum_{i=1}^{n}p_{i}\left\vert a_{i}\right\vert ^{2}-r^{2}>0,
\end{equation*}%
we may divide (\ref{5.3}) by $\ \sqrt{\sum_{i=1}^{n}p_{i}\left\vert
a_{i}\right\vert ^{2}-r^{2}}>0,$ getting%
\begin{equation}
\frac{\sum_{i=1}^{n}p_{i}\left\vert b_{i}\right\vert ^{2}}{\sqrt{%
\sum_{i=1}^{n}p_{i}\left\vert a_{i}\right\vert ^{2}-r^{2}}}+\sqrt{%
\sum_{i=1}^{n}p_{i}\left\vert a_{i}\right\vert ^{2}-r^{2}}\leq \frac{2\func{%
Re}\left( \sum_{i=1}^{n}p_{i}a_{i}b_{i}\right) }{\sqrt{\sum_{i=1}^{n}p_{i}%
\left\vert a_{i}\right\vert ^{2}-r^{2}}}.  \label{5.4}
\end{equation}%
Utilising the elementary inequality%
\begin{equation}
\frac{p}{\alpha }+q\alpha \geq 2\sqrt{pq}\text{ \ for \ }p,q\geq 0\text{ \
and \ }\alpha >0,  \label{5.4a}
\end{equation}%
we may write that%
\begin{equation}
2\sqrt{\sum_{i=1}^{n}p_{i}\left\vert b_{i}\right\vert ^{2}}\leq \frac{%
\sum_{i=1}^{n}p_{i}\left\vert b_{i}\right\vert ^{2}}{\sqrt{%
\sum_{i=1}^{n}p_{i}\left\vert a_{i}\right\vert ^{2}-r^{2}}}+\sqrt{%
\sum_{i=1}^{n}p_{i}\left\vert a_{i}\right\vert ^{2}-r^{2}}.  \label{5.5}
\end{equation}%
Combining (\ref{5.4a}) with (\ref{5.5}) we deduce%
\begin{equation}
\sqrt{\sum_{i=1}^{n}p_{i}\left\vert b_{i}\right\vert ^{2}}\leq \frac{\func{Re%
}\left( \sum_{i=1}^{n}p_{i}a_{i}b_{i}\right) }{\sqrt{\sum_{i=1}^{n}p_{i}%
\left\vert a_{i}\right\vert ^{2}-r^{2}}}.  \label{5.6}
\end{equation}%
Taking the square in (\ref{5.6}), we obtain%
\begin{equation*}
\sum_{i=1}^{n}p_{i}\left\vert b_{i}\right\vert ^{2}\left(
\sum_{i=1}^{n}p_{i}\left\vert a_{i}\right\vert ^{2}-r^{2}\right) \leq \left[ 
\func{Re}\left( \sum_{i=1}^{n}p_{i}a_{i}b_{i}\right) \right] ^{2},
\end{equation*}%
giving the third inequality in (\ref{5.2}).

The other inequalities are obvious.

To prove the sharpness of the constant, assume, under the hypothesis of the
theorem, that there exists a constant $c>0$ such that%
\begin{equation}
\sum_{i=1}^{n}p_{i}\left\vert a_{i}\right\vert
^{2}\sum_{i=1}^{n}p_{i}\left\vert b_{i}\right\vert ^{2}-\left[ \func{Re}%
\left( \sum_{i=1}^{n}p_{i}a_{i}b_{i}\right) \right] ^{2}\leq
cr^{2}\sum_{i=1}^{n}p_{i}\left\vert b_{i}\right\vert ^{2},  \label{5.7}
\end{equation}%
provided%
\begin{equation*}
\sum_{i=1}^{n}p_{i}\left\vert b_{i}-\overline{a_{i}}\right\vert ^{2}\leq
r^{2}<\sum_{i=1}^{n}p_{i}\left\vert a_{i}\right\vert ^{2}.
\end{equation*}%
Let $r=\sqrt{\varepsilon },$ $\varepsilon \in \left( 0,1\right) ,$ $%
a_{i},e_{i}\in \mathbb{C}$, $i\in \left\{ 1,\dots ,n\right\} $ with $%
\sum_{i=1}^{n}p_{i}\left\vert a_{i}\right\vert
^{2}=\sum_{i=1}^{n}p_{i}\left\vert e_{i}\right\vert ^{2}=1$ and $%
\sum_{i=1}^{n}p_{i}a_{i}e_{i}=0.$ Put $b_{i}=\overline{a_{i}}+\sqrt{%
\varepsilon }e_{i}.$ Then, obviously%
\begin{equation*}
\sum_{i=1}^{n}p_{i}\left\vert b_{i}-\overline{a_{i}}\right\vert ^{2}=r^{2},\
\ \ \ \sum_{i=1}^{n}p_{i}\left\vert a_{i}\right\vert ^{2}=1>r
\end{equation*}%
and 
\begin{gather*}
\sum_{i=1}^{n}p_{i}\left\vert b_{i}\right\vert
^{2}=\sum_{i=1}^{n}p_{i}\left\vert a_{i}\right\vert ^{2}+\varepsilon
\sum_{i=1}^{n}p_{i}\left\vert e_{i}\right\vert ^{2}=1+\varepsilon , \\
\func{Re}\left( \sum_{i=1}^{n}p_{i}a_{i}b_{i}\right)
=\sum_{i=1}^{n}p_{i}\left\vert a_{i}\right\vert ^{2}=1
\end{gather*}%
and thus%
\begin{equation*}
\sum_{i=1}^{n}p_{i}\left\vert a_{i}\right\vert
^{2}\sum_{i=1}^{n}p_{i}\left\vert b_{i}\right\vert ^{2}-\left[ \func{Re}%
\left( \sum_{i=1}^{n}p_{i}a_{i}b_{i}\right) \right] ^{2}=\varepsilon .
\end{equation*}%
Using (\ref{5.7}), we may write%
\begin{equation*}
\varepsilon \leq c\varepsilon \left( 1+\varepsilon \right) \text{ \ for \ }%
\varepsilon \in \left( 0,1\right) ,
\end{equation*}%
giving $1\leq c\left( 1+\varepsilon \right) $ for $\varepsilon \in \left(
0,1\right) .$ Making $\varepsilon \rightarrow 0+,$ we deduce $c\geq 1.$
\end{proof}

The following result also holds.

\begin{theorem}
\label{t5.2}Let $\mathbf{\bar{x}}=\left( x_{1},\dots ,x_{n}\right) $, $%
\mathbf{\bar{y}}=\left( y_{1},\dots ,y_{n}\right) \in \mathbb{K}^{n},$ $%
\mathbf{\bar{p}}=\left( p_{1},\dots ,p_{n}\right) \in \mathbb{R}_{+}^{n}$
with $\sum_{i=1}^{n}p_{i}=1$ and $\gamma ,\Gamma \in \mathbb{K}$ such that $%
\func{Re}\left( \gamma \bar{\Gamma}\right) >0$ and either%
\begin{equation}
\sum_{i=1}^{n}p_{i}\func{Re}\left[ \left( \Gamma \overline{y_{i}}%
-x_{i}\right) \left( \overline{x_{i}}-\bar{\gamma}y_{i}\right) \right] \geq
0,  \label{5.8}
\end{equation}%
or, equivalently,%
\begin{equation}
\sum_{i=1}^{n}p_{i}\left\vert x_{i}-\frac{\gamma +\Gamma }{2}\cdot \overline{%
y_{i}}\right\vert ^{2}\leq \frac{1}{4}\left\vert \Gamma -\gamma \right\vert
^{2}\sum_{i=1}^{n}p_{i}\left\vert y_{i}\right\vert ^{2}.  \label{5.9}
\end{equation}%
Then we have the inequalities%
\begin{align}
\sum_{i=1}^{n}p_{i}\left\vert x_{i}\right\vert
^{2}\sum_{i=1}^{n}p_{i}\left\vert y_{i}\right\vert ^{2}& \leq \frac{1}{4}%
\cdot \frac{\left\{ \func{Re}\left[ \left( \bar{\Gamma}+\bar{\gamma}\right)
\sum_{i=1}^{n}p_{i}x_{i}y_{i}\right] \right\} ^{2}}{\func{Re}\left( \Gamma 
\bar{\gamma}\right) }  \label{5.10} \\
& \leq \frac{1}{4}\cdot \frac{\left\vert \Gamma +\gamma \right\vert ^{2}}{%
\func{Re}\left( \Gamma \bar{\gamma}\right) }\left\vert
\sum_{i=1}^{n}p_{i}x_{i}y_{i}\right\vert ^{2}.  \notag
\end{align}%
The constant $\frac{1}{4}$ is best possible in both inequalities.
\end{theorem}

\begin{proof}
Define $b_{i}=x_{i}$ and $a_{i}=\frac{\bar{\Gamma}+\bar{\gamma}}{2}\cdot
y_{i}$ and $r=\frac{1}{2}\left\vert \Gamma -\gamma \right\vert \left(
\sum_{i=1}^{n}p_{i}\left\vert y_{i}\right\vert ^{2}\right) ^{\frac{1}{2}}.$
Then, by (\ref{5.9})%
\begin{align*}
\sum_{i=1}^{n}p_{i}\left\vert b_{i}-\overline{a_{i}}\right\vert ^{2}&
=\sum_{i=1}^{n}p_{i}\left\vert x_{i}-\frac{\gamma +\Gamma }{2}\cdot 
\overline{y_{i}}\right\vert ^{2} \\
& \leq \frac{1}{4}\left\vert \Gamma -\gamma \right\vert
^{2}\sum_{i=1}^{n}p_{i}\left\vert y_{i}\right\vert ^{2}=r^{2},
\end{align*}%
showing that the first condition in (\ref{5.1}) is fulfilled.

We also have%
\begin{align*}
\sum_{i=1}^{n}p_{i}\left\vert a_{i}\right\vert ^{2}-r^{2}&
=\sum_{i=1}^{n}p_{i}\left\vert \frac{\Gamma +\gamma }{2}\right\vert
^{2}\left\vert y_{i}\right\vert ^{2}-\frac{1}{4}\left\vert \Gamma -\gamma
\right\vert ^{2}\sum_{i=1}^{n}p_{i}\left\vert y_{i}\right\vert ^{2} \\
& =\func{Re}\left( \Gamma \bar{\gamma}\right) \sum_{i=1}^{n}p_{i}\left\vert
y_{i}\right\vert ^{2}>0
\end{align*}%
since $\func{Re}\left( \gamma \bar{\Gamma}\right) >0,$ and thus the
condition in (\ref{5.1}) is also satisfied.

Using the second inequality in (\ref{5.2}), one may write%
\begin{eqnarray*}
&&\sum_{i=1}^{n}p_{i}\left\vert \frac{\Gamma +\gamma }{2}\right\vert
^{2}\left\vert y_{i}\right\vert ^{2}\sum_{i=1}^{n}p_{i}\left\vert
x_{i}\right\vert ^{2}-\left[ \func{Re}\sum_{i=1}^{n}p_{i}\left( \frac{\bar{%
\Gamma}+\bar{\gamma}}{2}\right) y_{i}x_{i}\right] ^{2} \\
&\leq &\frac{1}{4}\left\vert \Gamma -\gamma \right\vert
^{2}\sum_{i=1}^{n}p_{i}\left\vert y_{i}\right\vert
^{2}\sum_{i=1}^{n}p_{i}\left\vert x_{i}\right\vert ^{2},
\end{eqnarray*}%
giving%
\begin{equation*}
\frac{\left\vert \Gamma +\gamma \right\vert ^{2}-\left\vert \Gamma -\gamma
\right\vert ^{2}}{4}\sum_{i=1}^{n}p_{i}\left\vert y_{i}\right\vert
^{2}\sum_{i=1}^{n}p_{i}\left\vert x_{i}\right\vert ^{2}\leq \frac{1}{4}\func{%
Re}\left[ \left( \bar{\Gamma}+\bar{\gamma}\right)
\sum_{i=1}^{n}p_{i}x_{i}y_{i}\right] ^{2},
\end{equation*}%
which is clearly equivalent to the first inequality in (\ref{5.10}).

The second inequality in (\ref{5.10}) is obvious.

To prove the sharpness of the constant $\frac{1}{4},$ assume that the first
inequality in (\ref{5.10}) holds with a constant $c>0,$ i.e.,%
\begin{equation}
\sum_{i=1}^{n}p_{i}\left\vert x_{i}\right\vert
^{2}\sum_{i=1}^{n}p_{i}\left\vert y_{i}\right\vert ^{2}\leq C\cdot \frac{%
\left\{ \func{Re}\left[ \left( \bar{\Gamma}+\bar{\gamma}\right)
\sum_{i=1}^{n}p_{i}x_{i}y_{i}\right] \right\} ^{2}}{\func{Re}\left( \Gamma 
\bar{\gamma}\right) },  \label{5.11}
\end{equation}%
provided $\func{Re}\left( \gamma \bar{\Gamma}\right) >0$ and either (\ref%
{5.8}) or (\ref{5.9}) holds.

Assume that $\Gamma ,\gamma >0$ and let $x_{i}=\gamma \bar{y}_{i}.$ Then (%
\ref{5.8}) holds true and by (\ref{5.11}) we deduce%
\begin{equation*}
\gamma ^{2}\left( \sum_{i=1}^{n}p_{i}\left\vert y_{i}\right\vert ^{2}\right)
^{2}\leq C\frac{\left( \Gamma +\gamma \right) ^{2}\gamma ^{2}\left(
\sum_{i=1}^{n}p_{i}\left\vert y_{i}\right\vert ^{2}\right) ^{2}}{\Gamma
\gamma },
\end{equation*}%
giving%
\begin{equation}
\Gamma \gamma \leq C\left( \Gamma +\gamma \right) ^{2}\text{ \ for any \ }%
\Gamma ,\gamma >0.  \label{5.12}
\end{equation}%
Let $\varepsilon \in \left( 0,1\right) $ and choose in (\ref{5.12}) $\Gamma
=1+\varepsilon ,$ $\gamma =1-\varepsilon >0$ to get $1-\varepsilon ^{2}\leq
4C$ for any $\varepsilon \in \left( 0,1\right) .$ Letting $\varepsilon
\rightarrow 0+,$ we deduce $C\geq \frac{1}{4}$ and the sharpness of the
constant is proved.

Finally, we note that the conditions (\ref{5.8}) and (\ref{5.9}) are
equivalent since in an inner product space $\left( H,\left\langle \cdot
,\cdot \right\rangle \right) $ for any vectors $x,z,Z\in H$ one has $\func{Re%
}\left\langle Z-x,x-z\right\rangle \geq 0$ iff $\left\Vert x-\frac{z+Z}{2}%
\right\Vert \leq \frac{1}{2}\left\Vert Z-z\right\Vert $ \cite{1b}. We omit
the details.
\end{proof}

\section{More Reverses of the (CBS)-Inequality}

The following result holds.

\begin{theorem}
\label{t6.1}Let $\mathbf{\bar{a}}=\left( a_{1},\dots ,a_{n}\right) $, $%
\mathbf{\bar{b}}=\left( b_{1},\dots ,b_{n}\right) \in \mathbb{K}^{n}$ and $%
\mathbf{\bar{p}}=\left( p_{1},\dots ,p_{n}\right) \in \mathbb{R}_{+}^{n}$
with $\sum_{i=1}^{n}p_{i}=1.$ If $r>0$ and the following condition is
satisfied%
\begin{equation}
\sum_{i=1}^{n}p_{i}\left\vert b_{i}-\overline{a_{i}}\right\vert ^{2}\leq
r^{2},  \label{6.1}
\end{equation}%
then we have the inequalities%
\begin{align}
0& \leq \left( \sum_{i=1}^{n}p_{i}\left\vert b_{i}\right\vert
^{2}\sum_{i=1}^{n}p_{i}\left\vert a_{i}\right\vert ^{2}\right) ^{\frac{1}{2}%
}-\left\vert \sum_{i=1}^{n}p_{i}a_{i}b_{i}\right\vert  \label{6.2} \\
& \leq \left( \sum_{i=1}^{n}p_{i}\left\vert b_{i}\right\vert
^{2}\sum_{i=1}^{n}p_{i}\left\vert a_{i}\right\vert ^{2}\right) ^{\frac{1}{2}%
}-\left\vert \sum_{i=1}^{n}p_{i}\func{Re}\left( a_{i}b_{i}\right) \right\vert
\notag \\
& \leq \left( \sum_{i=1}^{n}p_{i}\left\vert b_{i}\right\vert
^{2}\sum_{i=1}^{n}p_{i}\left\vert a_{i}\right\vert ^{2}\right) ^{\frac{1}{2}%
}-\sum_{i=1}^{n}p_{i}\func{Re}\left( a_{i}b_{i}\right)  \notag \\
& \leq \frac{1}{2}r^{2}.  \notag
\end{align}%
The constant $\frac{1}{2}$ is best possible in (\ref{6.2}) in the sense that
it cannot be replaced by a smaller constant.
\end{theorem}

\begin{proof}
The condition (\ref{6.1}) is clearly equivalent to%
\begin{equation}
\sum_{i=1}^{n}p_{i}\left\vert b_{i}\right\vert
^{2}+\sum_{i=1}^{n}p_{i}\left\vert a_{i}\right\vert ^{2}\leq
2\sum_{i=1}^{n}p_{i}\func{Re}\left( b_{i}a_{i}\right) +r^{2}.  \label{6.3}
\end{equation}%
Using the elementary inequality%
\begin{equation}
2\left( \sum_{i=1}^{n}p_{i}\left\vert b_{i}\right\vert
^{2}\sum_{i=1}^{n}p_{i}\left\vert a_{i}\right\vert ^{2}\right) ^{\frac{1}{2}%
}\leq \sum_{i=1}^{n}p_{i}\left\vert b_{i}\right\vert
^{2}+\sum_{i=1}^{n}p_{i}\left\vert a_{i}\right\vert ^{2}  \label{6.4}
\end{equation}%
and (\ref{6.3}), we deduce%
\begin{equation}
2\left( \sum_{i=1}^{n}p_{i}\left\vert b_{i}\right\vert
^{2}\sum_{i=1}^{n}p_{i}\left\vert a_{i}\right\vert ^{2}\right) ^{\frac{1}{2}%
}\leq 2\sum_{i=1}^{n}p_{i}\func{Re}\left( b_{i}a_{i}\right) +r^{2},
\label{6.5}
\end{equation}%
giving the last inequality in (\ref{6.2}). The other inequalities are
obvious.

To prove the sharpness of the constant $\frac{1}{2},$ assume that%
\begin{equation}
0\leq \left( \sum_{i=1}^{n}p_{i}\left\vert b_{i}\right\vert
^{2}\sum_{i=1}^{n}p_{i}\left\vert a_{i}\right\vert ^{2}\right) ^{\frac{1}{2}%
}-\sum_{i=1}^{n}p_{i}\func{Re}\left( b_{i}a_{i}\right) \leq cr^{2}
\label{6.6}
\end{equation}%
for any $\mathbf{\bar{a}}$, $\mathbf{\bar{b}}\in \mathbb{K}^{n}$ and $r>0$
satisfying (\ref{6.1}).

Assume that $\mathbf{\bar{a}}$, $\mathbf{\bar{e}}\in H,$ $\mathbf{\bar{e}}%
=\left( e_{1},\dots ,e_{n}\right) $ with $\sum_{i=1}^{n}p_{i}\left\vert
a_{i}\right\vert ^{2}=\sum_{i=1}^{n}p_{i}\left\vert e_{i}\right\vert ^{2}=1$
and $\sum_{i=1}^{n}p_{i}a_{i}e_{i}=0.$ If $r=\sqrt{\varepsilon },$ $%
\varepsilon >0,$ and if we define $\mathbf{\bar{b}}=\overline{\mathbf{\bar{a}%
}}+\sqrt{\varepsilon }\mathbf{\bar{e}}$ where $\overline{\mathbf{\bar{a}}}%
=\left( \overline{a_{1}},\dots ,\overline{a_{n}}\right) \in \mathbb{K}^{n},$
then $\sum_{i=1}^{n}p_{i}\left\vert b_{i}-\overline{a_{i}}\right\vert
^{2}=\varepsilon =r^{2},$ showing that the condition (\ref{6.1}) is
fulfilled.

On the other hand,%
\begin{align*}
& \left( \sum_{i=1}^{n}p_{i}\left\vert b_{i}\right\vert
^{2}\sum_{i=1}^{n}p_{i}\left\vert a_{i}\right\vert ^{2}\right) ^{\frac{1}{2}%
}-\sum_{i=1}^{n}p_{i}\func{Re}\left( b_{i}a_{i}\right)  \\
& =\left( \sum_{i=1}^{n}p_{i}\left\vert \overline{a_{i}}+\sqrt{\varepsilon }%
e_{i}\right\vert ^{2}\right) ^{\frac{1}{2}}-\sum_{i=1}^{n}p_{i}\func{Re}%
\left[ \left( \overline{a_{i}}+\sqrt{\varepsilon }e_{i}\right) a_{i}\right] 
\\
& =\left( \sum_{i=1}^{n}p_{i}\left\vert a_{i}\right\vert ^{2}+\varepsilon
\sum_{i=1}^{n}\left\vert e_{i}\right\vert ^{2}\right) ^{\frac{1}{2}%
}-\sum_{i=1}^{n}p_{i}\left\vert a_{i}\right\vert ^{2} \\
& =\sqrt{1+\varepsilon }-1.
\end{align*}%
Utilizing (\ref{6.6}), we conclude that%
\begin{equation}
\sqrt{1+\varepsilon }-1\leq c\varepsilon \text{ \ for any \ }\varepsilon >0.
\label{6.7}
\end{equation}%
Multiplying (\ref{6.7}) by $\sqrt{1+\varepsilon }+1>0$ and thus dividing by $%
\varepsilon >0,$ we get%
\begin{equation}
\left( \sqrt{1+\varepsilon }-1\right) c\geq 1\text{ \ for any \ }\varepsilon
>0.  \label{6.8}
\end{equation}%
Letting $\varepsilon \rightarrow 0+$ in (\ref{6.8}), we deduce $c\geq \frac{1%
}{2},$ and the theorem is proved.
\end{proof}

Finally, the following result also holds.

\begin{theorem}
\label{t6.2}Let $\mathbf{\bar{x}}=\left( x_{1},\dots ,x_{n}\right) $, $%
\mathbf{\bar{y}}=\left( y_{1},\dots ,y_{n}\right) \in \mathbb{K}^{n},$ $%
\mathbf{\bar{p}}=\left( p_{1},\dots ,p_{n}\right) \in \mathbb{R}_{+}^{n}$
with $\sum_{i=1}^{n}p_{i}=1,$ and $\gamma ,\Gamma \in \mathbb{K}$ with $%
\Gamma \neq \gamma ,-\gamma ,$ so that either%
\begin{equation}
\sum_{i=1}^{n}p_{i}\func{Re}\left[ \left( \Gamma \overline{y_{i}}%
-x_{i}\right) \left( \overline{x_{i}}-\bar{\gamma}y_{i}\right) \right] \geq
0,  \label{6.9}
\end{equation}%
or, equivalently,%
\begin{equation}
\sum_{i=1}^{n}p_{i}\left\vert x_{i}-\frac{\gamma +\Gamma }{2}\cdot \overline{%
y_{i}}\right\vert ^{2}\leq \frac{1}{4}\left\vert \Gamma -\gamma \right\vert
^{2}\sum_{i=1}^{n}p_{i}\left\vert y_{i}\right\vert ^{2}  \label{6.10}
\end{equation}%
holds. Then we have the inequalities%
\begin{align}
0& \leq \left( \sum_{i=1}^{n}p_{i}\left\vert x_{i}\right\vert
^{2}\sum_{i=1}^{n}p_{i}\left\vert y_{i}\right\vert ^{2}\right) ^{\frac{1}{2}%
}-\left\vert \sum_{i=1}^{n}p_{i}x_{i}y_{i}\right\vert  \label{6.11} \\
& \leq \left( \sum_{i=1}^{n}p_{i}\left\vert x_{i}\right\vert
^{2}\sum_{i=1}^{n}p_{i}\left\vert y_{i}\right\vert ^{2}\right) ^{\frac{1}{2}%
}-\left\vert \sum_{i=1}^{n}p_{i}\func{Re}\left[ \frac{\bar{\Gamma}+\bar{%
\gamma}}{\left\vert \Gamma +\gamma \right\vert }x_{i}y_{i}\right] \right\vert
\notag \\
& \leq \left( \sum_{i=1}^{n}p_{i}\left\vert x_{i}\right\vert
^{2}\sum_{i=1}^{n}p_{i}\left\vert y_{i}\right\vert ^{2}\right) ^{\frac{1}{2}%
}-\sum_{i=1}^{n}p_{i}\func{Re}\left[ \frac{\bar{\Gamma}+\bar{\gamma}}{%
\left\vert \Gamma +\gamma \right\vert }x_{i}y_{i}\right]  \notag \\
& \leq \frac{1}{4}\cdot \frac{\left\vert \Gamma -\gamma \right\vert ^{2}}{%
\left\vert \Gamma +\gamma \right\vert }\sum_{i=1}^{n}p_{i}\left\vert
y_{i}\right\vert ^{2}.  \notag
\end{align}%
The constant $\frac{1}{4}$ in the last inequality is best possible.
\end{theorem}

\begin{proof}
Consider $b_{i}=x_{i},$ $a_{i}=\frac{\bar{\Gamma}+\bar{\gamma}}{2}\cdot
y_{i},$ $i\in \left\{ 1,\dots ,n\right\} $ and 
\begin{equation*}
r:=\frac{1}{2}\left( \Gamma -\gamma \right) \left(
\sum_{i=1}^{n}p_{i}\left\vert y_{i}\right\vert ^{2}\right) ^{\frac{1}{2}}.
\end{equation*}%
Then, by (\ref{6.10}), we have%
\begin{equation*}
\sum_{i=1}^{n}p_{i}\left\vert b_{i}-\overline{a_{i}}\right\vert
^{2}=\sum_{i=1}^{n}p_{i}\left\vert x_{i}-\frac{\gamma +\Gamma }{2}\cdot
y_{i}\right\vert ^{2}\leq \frac{1}{4}\left\vert \Gamma -\gamma \right\vert
^{2}\sum_{i=1}^{n}p_{i}\left\vert y_{i}\right\vert ^{2}=r^{2}
\end{equation*}%
showing that (\ref{6.1}) is valid.

By the use of the last inequality in (\ref{6.2}), we have%
\begin{align*}
0& \leq \left( \sum_{i=1}^{n}p_{i}\left\vert x_{i}\right\vert
^{2}\sum_{i=1}^{n}p_{i}\left\vert \frac{\Gamma +\gamma }{2}\right\vert
^{2}\left\vert y_{i}\right\vert ^{2}\right) ^{\frac{1}{2}}-%
\sum_{i=1}^{n}p_{i}\func{Re}\left[ \frac{\bar{\Gamma}+\bar{\gamma}}{2}%
x_{i}y_{i}\right] \\
& \leq \frac{1}{8}\left\vert \Gamma -\gamma \right\vert
^{2}\sum_{i=1}^{n}p_{i}\left\vert y_{i}\right\vert ^{2}.
\end{align*}%
Dividing by $\frac{1}{2}\left\vert \Gamma +\gamma \right\vert >0,$ we deduce%
\begin{align*}
0& \leq \left( \sum_{i=1}^{n}p_{i}\left\vert x_{i}\right\vert
^{2}\sum_{i=1}^{n}p_{i}\left\vert y_{i}\right\vert ^{2}\right) ^{\frac{1}{2}%
}-\sum_{i=1}^{n}p_{i}\func{Re}\left[ \frac{\bar{\Gamma}+\bar{\gamma}}{%
\left\vert \Gamma +\gamma \right\vert }x_{i}y_{i}\right] \\
& \leq \frac{1}{4}\cdot \frac{\left\vert \Gamma -\gamma \right\vert ^{2}}{%
\left\vert \Gamma +\gamma \right\vert }\sum_{i=1}^{n}p_{i}\left\vert
y_{i}\right\vert ^{2},
\end{align*}%
which is the last inequality in (\ref{6.11}).

The other inequalities are obvious.

To prove the sharpness of the constant $\frac{1}{4},$ assume that there
exists a constant $c>0,$ such that%
\begin{eqnarray}
&&\left( \sum_{i=1}^{n}p_{i}\left\vert x_{i}\right\vert
^{2}\sum_{i=1}^{n}p_{i}\left\vert y_{i}\right\vert ^{2}\right) ^{\frac{1}{2}%
}-\sum_{i=1}^{n}p_{i}\func{Re}\left[ \frac{\bar{\Gamma}+\bar{\gamma}}{%
\left\vert \Gamma +\gamma \right\vert }x_{i}y_{i}\right]   \label{6.12} \\
&\leq &c\cdot \frac{\left\vert \Gamma -\gamma \right\vert ^{2}}{\left\vert
\Gamma +\gamma \right\vert }\sum_{i=1}^{n}p_{i}\left\vert y_{i}\right\vert
^{2},  \notag
\end{eqnarray}%
provided either (\ref{6.9}) or (\ref{6.10}) holds.

Let $n=2,$ $\mathbf{\bar{y}}=\left( 1,1\right) ,$ $\mathbf{\bar{x}}=\left(
x_{1},x_{2}\right) \in \mathbb{R}^{2},$ $\mathbf{\bar{p}}=\left( \frac{1}{2},%
\frac{1}{2}\right) $ and $\Gamma ,\gamma >0$ with $\Gamma >\gamma .$ Then by
(\ref{6.12}) we deduce%
\begin{equation}
\sqrt{2}\sqrt{x_{1}^{2}+x_{2}^{2}}-\left( x_{1}+x_{2}\right) \leq 2c\frac{%
\left( \Gamma -\gamma \right) ^{2}}{\Gamma +\gamma }.  \label{6.13}
\end{equation}%
If $x_{1}=\Gamma ,$ $x_{2}=\gamma ,$ then $\left( \Gamma -x_{1}\right)
\left( x_{1}-\gamma \right) +\left( \Gamma -x_{2}\right) \left( x_{2}-\gamma
\right) =0,$ showing that the condition (\ref{6.9}) is valid for $n=2$ and $%
\mathbf{\bar{p}}$, $\mathbf{\bar{x}}$, $\mathbf{\bar{y}}$ as above.
Replacing $x_{1}$ and $x_{2}$ in (\ref{6.13}), we deduce%
\begin{equation}
\sqrt{2}\sqrt{\Gamma ^{2}+\gamma ^{2}}-\left( \Gamma +\gamma \right) \leq 2c%
\frac{\left( \Gamma -\gamma \right) ^{2}}{\Gamma +\gamma }.  \label{6.14}
\end{equation}%
If in (\ref{6.14}) we choose $\Gamma =1+\varepsilon ,$ $\gamma
=1-\varepsilon $ with $\varepsilon \in \left( 0,1\right) ,$ we deduce%
\begin{equation}
\sqrt{1+\varepsilon ^{2}}-1\leq 2c\varepsilon ^{2}.  \label{6.15}
\end{equation}%
Finally, multiplying (\ref{6.15}) with $\sqrt{1+\varepsilon ^{2}}+1>0$ and
then dividing by $\varepsilon ^{2},$ we deduce%
\begin{equation}
1\leq 2c\left( \sqrt{1+\varepsilon ^{2}}+1\right) \text{ \ for any \ }%
\varepsilon >0.  \label{6.16}
\end{equation}%
Letting $\varepsilon \rightarrow 0+$ in (\ref{6.16}), we get $c\geq \frac{1}{%
4},$ and the sharpness of the constant is proved.
\end{proof}

\begin{remark}
The integral version may be stated in a canonical way. The corresponding
inequalities for integrals will be considered in another  work devoted to
positive linear functionals with complex values \ that is in preparation.
\end{remark}

\end{document}